\documentclass[12pt,a4paper,reqno]{amsart}  

\usepackage{mathtools}
\usepackage{url}

\usepackage[headings]{fullpage}
\usepackage{comment} 
\usepackage{datetime}  
\usepackage[english]{babel} 

\usepackage[full]{textcomp}
\usepackage[osf]{newtxtext} 
\usepackage{cabin} 
\usepackage[bigdelims,vvarbb]{newtxmath} 
\usepackage[cal=boondoxo]{mathalfa} 
\usepackage{zlmtt}

\usepackage{microtype} 
 
\allowdisplaybreaks

\renewcommand{\qedsymbol}{$\square$}

\newtheoremstyle{sltheorems}
{10pt}
{6pt}
{\slshape}
{}
{\bfseries}
{.}
{.5em}
{\thmname{#1}\thmnumber{ #2}\thmnote{ (#3)}}

\theoremstyle{sltheorems} 
\newtheorem{teor}{Theorem}
\newtheorem{lem}{Lemma} 
\newtheorem{cor}{Corollary} 

\newtheoremstyle{remark}
{10pt}
{6pt}
{\rm} 
{}
{\bfseries}
{.}
{.5em}
{\thmname{#1}\thmnumber{ #2}\thmnote{ (#3)}}
 \theoremstyle{remark}

\makeatletter

\def\env@Biggcases{%
  \let\@ifnextchar\new@ifnextchar
  \Biggl\lbrace
  \def\arraystretch{1.2}%
  \array{@{}l@{\quad}l@{}}%
}
\makeatother

\newcommand{\z}{\mathbb}

\newcommand{\dd}{\mathrm d}

\begin{document}


\baselineskip=17pt


\title[Diophantine approximation with primes]{Diophantine approximation with one prime, two squares of primes and one $k$-th power of a prime}

\author{Alessandro Gambini
}

\date{}

\begin{abstract}
Let $1<k<14/5$, $\lambda_1,\lambda_2,\lambda_3$ and $\lambda_4$ be non-zero real numbers, not all of the same sign such that $\lambda_1/\lambda_2$ is irrational and let $\omega$ be a real number. We prove that the inequality $|\lambda_1p_1+\lambda_2p_2^2+\lambda_3p_3^2+\lambda_4p_4^k-\omega|\le (\max_j p_j)^{-\frac{14-5k}{28k}+\varepsilon}$ has infinitely many solutions in prime variables $p_1,p_2,p_3,p_4$ for any $\varepsilon>0$.
\end{abstract}

\maketitle


\renewcommand{\thefootnote}{}



\renewcommand{\thefootnote}{\arabic{footnote}}
\setcounter{footnote}{0}


\section{Introduction}
This paper deals with a Diophantine inequality with prime variables involving a prime, two squares of primes and one $k$-th power of a prime. In particular we prove the following theorem:

\begin{teor}
Assume that $1<k<14/5$, $\lambda_1,\lambda_2,\lambda_3$ and $\lambda_4$ be non-zero real numbers, not all of the same sign, that $\lambda_1/\lambda_2$ is irrational and let $\omega$ be a real number. The inequality
\begin{align}\label{teorema}
\left|\lambda_1p_1+\lambda_2p_2^2+\lambda_3p_3^2+\lambda_4p_4^k-\omega\right|\le\left(\max_j p_j\right)^{-\frac{14-5k}{28k}+\varepsilon}
\end{align}
has infinitely many solutions in prime variables $p_1,p_2,p_3,p_4$ for any $\varepsilon>0$.
\end{teor}

Many recent such results are known with various types of assumptions and conclusions. Many of them deal with the number of \emph{exceptional} real numbers $\omega$ such that the inequality
\begin{align*}
|\lambda_1p_1^{k_1}+\cdots+\lambda_rp_r^{k_r}-\omega|\le\eta
\end{align*} 
has no solution in prime variables $p_1,\ldots,p_r$, for small $\eta>0$ fixed.

Br\"udern, Cook and Perelli in \cite{brudern-cook-perelli} dealt with binary linear forms in prime arguments, Cook and Fox in \cite{cook-fox} dealt with a ternary form with squares of primes that was improved in term of approximation by Harman in \cite{harman2004ternary}. Cook in \cite{cook2001general} gave a more general description of the problem, later improved by Cook and Harman in \cite{cook-harman}.

There are some differences between the results quoted above and our purpose: in our case the value of $\eta$ does depend on the primes $p_j$ and it will be actually a negative power of the maximum of the $p_j$ while in the papers quoted above $\eta$ is a small negative power of $\omega$. In their papers the assumption that the coefficients $\lambda_j$ are all positive is not a restriction. Moreover $k_j$ is the same positive integer for all $j$. Nevertheless the assumption that $\lambda_1/\lambda_2$ must be irrational is still the heart of the matter. 

Vaughan in \cite{vaughan1974diophantineI} follows another approach, that is the same we are using in our article: dealing with a ternary linear form in prime arguments and assuming some more suitable conditions on the $\lambda_j$, he proved that there are infinitely many solutions of the problem 
\begin{align*}
|\lambda_1p_1^{k_1}+\cdots+\lambda_rp_r^{k_r}-\omega|\le\eta
\end{align*} 
when $\eta$ depends on the maximum of the $p_j$; in his case $\eta=(\max_j p_j)^{-\frac{1}{10}}$. Such result was improved by Baker and Harman in \cite{baker-harman} with exponent $-\frac16$, by Harman in \cite{harman1991general} with exponent $-\frac15$ and finally by Matom\"aki in \cite{matomaki} with exponent $-\frac29$.

Languasco and Zaccagnini in \cite{Languasco-Zaccagnini2016} and \cite{languasco-zaccagnini-ternary} dealt with a ternary problem with a $k$-th power of a prime. In this case the value of $\eta$ is a negative power of the maximum of the $p_j$ also depending on the parameter $k$: the idea in this case is to get both the widest $k$-range and the strongest bound for the approximation.

Languasco and Zaccagnini dealt also with a quaternary form in \cite{languasco-zaccagnini-quaternary} with a prime and 3 squares of primes obtaining $\eta=(\max_j p_j)^{-\frac{1}{18}}$; this was improved by Liu \& Sun in \cite{liu2013diophantine} with $\eta=(\max_j p_j)^{-\frac{1}{16}}$ using the Harman technique. Wang \& Yao in [WY] improved the approximation to the exponent $-\frac{1}{14}$; in this paper we generalized the problem to a real power $k\in\left(1,\frac{14}{5}\right)$.

\section{Outline of the proof}

We use a variant of the classical circle method that was introduced by Davenport and Heilbronn in 1946 \cite{davenport1946indefinite} in order to attack this kind of Diophantine problems. 
The integration on a circle, or equivalently on the interval $[0,1]$, is replaced by integration on the whole
real line. 

Throughout this paper $p_i$ denotes a prime number, $k\ge1$ is a real number, $\varepsilon$ is an arbitrarily small positive number whose value could vary depending on the occurrences and $\omega$ is a fixed real number. In order to prove that \eqref{teorema} has infinitely many solutions, it is sufficient to construct an increasing sequence $X_n$ that tends to infinity such that \eqref{teorema} has at least one solution with $\max p_j\in[\delta X_n,X_n]$, with $\delta>0$ fixed, depending on the choice of $\lambda_j$. Let $q$ be a denominator of a convergent to $\lambda_1/\lambda_2$ and let $X_n=X$ (dropping the suffix $n$) run through the sequence $X=q^{7/3}$.
Set
\begin{align}\label{S_k}
&S_k(\alpha)=\sum_{\delta X\le p^k\le X}\log p\  e(p^k\alpha)\\
&U_k(\alpha)=\sum_{\delta X\le p^k\le X}e(p^k\alpha)\label{U_k} \\
&T_k(\alpha)=\int_{(\delta X)^{\frac1k}}^{X^{\frac1k}}e(\alpha t^k)\,\dd t,\label{T_k}
\end{align}
where $e(\alpha)=e^{2\pi i\alpha}$. We will approximate $S_k$ with $T_k$ and $U_k$.

By the Prime Number Theorem and first derivative estimates for trigonometric integrals we have
\begin{align}\label{stima_tk}
S_k(\alpha)\ll X^{\frac1k}, \qquad T_k(\alpha)\ll_{k,\delta} X^{\frac1k-1}\min(X,|\alpha|^{-1}).
\end{align}

Moreover the Euler summation formula implies that
\begin{align}\label{t-u}
T_k(\alpha)-U_k(\alpha)\ll 1+|\alpha|X.
\end{align}

We also need a continuous function we will use to detect the solutions of \eqref{teorema}, so we introduce
\begin{align*}
\widehat{K}_\eta(\alpha):=\max\{0,\eta-|\alpha|\}\quad\mbox{where}\quad\eta>0
\end{align*}
whose inverse Fourier transform is
\begin{align*}
K_{\eta}(\alpha)=\left(\frac{\sin(\pi\alpha\eta)}{\pi\alpha}\right)^2
\end{align*}
for $\alpha\neq0$ and, by continuity, $K_{\eta}(0)=\eta^2$. It vanishes at infinity like $|\alpha|^{-2}$ and in fact it is trivial to prove that
\begin{align}\label{k_eta}
K_{\eta}(\alpha)\ll\min(\eta^2,|\alpha|^{-2}).
\end{align}

The original works of Davenport-Heillbronn in \cite{davenport1946indefinite} and later Vaughan in \cite{vaughan1974diophantineI} and \cite{vaughan1974diophantineII} approximate directly the difference $|S_k(\alpha)-T_k(\alpha)|$, estimating it with $O(1)$ using the Euler summation formula. The $L^2$-norm estimation approach (see \cite{brudern-cook-perelli} and \cite{Languasco-Zaccagnini2016}) improves these estimation taking the $L^2$-norm of $|S_k(\alpha)-T_k(\alpha)|$ leading to significantly better conditions and to have a wider major arc compared to the original DH approach. In fact, setting the generalized version of the Selberg integral
\begin{align*}
\mathcal{J}_k(X,h)=\int_{X}^{2X}\left(\theta((x+h)^{\frac1k})-\theta(x^{\frac1k})-((x+h)^{\frac1k}-x^{\frac1k})\right)^2\dd x,
\end{align*}
we have the following lemmas.

\begin{lem}[\cite{Languasco-Zaccagnini2016}, Theorem 1]\label{s-u}
Let $k\ge1$ be a real number. For $0<Y<\frac12$ we have
\begin{align*}
\int_{-Y}^{Y}|S_k(\alpha)-U_k(\alpha)|^2 \dd\alpha\ll_k\frac{{X^{\frac2k-2}}\log^2 X}{Y}+Y^2X+Y^2\mathcal{J}_k\left(X,\frac{1}{2Y}\right).
\end{align*}
\end{lem}
\begin{lem}[\cite{Languasco-Zaccagnini2016}, Theorem 2]\label{jk_stima}
Let $k\ge1$ be a real number and $\varepsilon$ be an arbitrarily small positive constant. There exists a positive constant $c_1(\varepsilon)$, which does not depend on $k$, such that
\begin{align*}
\mathcal{J}_k(X,h)\ll_k h^2X^{\frac2k-1}\exp\left(-c_1\left(\frac{\log X}{\log \log X}\right)^{\frac13}\right)
\end{align*}
uniformly for $X^{1-\frac{5}{6k}+\varepsilon}\le h\le X$.
\end{lem}

\subsection{Setting the problem}
Let 
\begin{align*}\mathcal{P}(X)=\{(p_1,p_2,p_3,p_4): \delta X<p_1<X,\ \delta X<p_2^2,p_3^2<X,\ \delta X<p_4^k<X\}
\end{align*} and let us define
\begin{align*}
\mathcal{I}(\eta,\omega,\mathfrak{X})=\int_{\mathfrak{X}}S_1(\lambda_1\alpha)S_2(\lambda_2\alpha)S_2(\lambda_3\alpha)S_k(\lambda_4\alpha)K_{\eta}(\alpha)e(-\omega\alpha)\dd\alpha
\end{align*}
where $\mathfrak{X}$ is a measurable subset of $\z{R}$.

From the definitions of the $S_j(\lambda_i\alpha)$ and performing the Fourier transform for $K_{\eta}(\alpha)$, we get
\begin{align*}
\mathcal{I}(\eta,\omega,\z{R})&=\sum_{\mathcal{P}(X)}\log p_1\log p_2\log p_3\log p_4\cdot \\&\quad\left(\max(0,\eta-\left|\lambda_1p_1+\lambda_2p_2^2+\lambda_3p_3^2+\lambda_4p_4^k-\omega\right|)\right)\nonumber\\ &\leq \eta(\log X)^4\mathcal{N}(X),
\end{align*}
where $\mathcal{N}(X)$ actually denotes the number of solutions of the inequality \eqref{teorema} with $(p_1,p_2,p_3,p_4)\in\mathcal{P}(X)$. In other words $\mathcal{I}(\eta,\omega,\z{R})$ provides a lower bound for the quantity we are interested in; therefore it is sufficient to prove that $\mathcal{I}(\eta,\omega,\z{R})>0$.

We now decompose $\z{R}$ into subsets such that $\z{R}=\mathcal{M}\,\cup\,m\,\cup\,t$ where $\mathcal{M}$ is the major arc, $m$ is the minor arc (or intermediate arc) and $t$ is the trivial arc. The decomposition is the following:

\begin{align*}
\mathcal{M}=\left[-\frac{P}{X},\frac{P}{X}\right] \qquad m=\left[\frac{P}{X},R\right]\cup\left[-R,-\frac{P}{X}\right]\qquad t=\z{R}\backslash(\mathcal{M}\cup m),
\end{align*}
so that  $\mathcal{I}(\eta,\omega,\z{R})=\mathcal{I}(\eta,\omega,\mathcal{M})+\mathcal{I}(\eta,\omega,m)+\mathcal{I}(\eta,\omega,t)$.

The parameters $P=P(X)>1$ and $R=R(X)>1/\eta$ are chosen later (see \eqref{cond_major} and \eqref{cond_trivial}) as well as $\eta=\eta(X)$, that, as we explained before, we would like to be a small negative power of $\max p_j$ (and so of $X$, see \eqref{final}).

We are expecting to have on $\mathcal{M}$ the main term with the right order of magnitude without any special hypothesis on the coefficients $\lambda_j$. It is necessary to prove that $\mathcal{I}(\eta,\omega,m)$ and $\mathcal{I}(\eta,\omega,t)$ are both $o(\mathcal{I}(\eta,\omega,\mathcal{M}))$: the contribution from the trivial is ``tiny" with respect to the main term. The real problem is on the minor arc where we will need the full force of the hypothesis on the $\lambda_j$ and the theory of continued fractions.

\textbf{Remark}: from now on, anytime we use the symbol $\ll$ or $\gg$ we drop the dependence of the approximation from the constants $\lambda_j, \delta$ and $k$.

\subsection{Lemmas}

In this paper we will also use Lemmas 3-4-10 of \cite{GLZ} that allow us to have an estimation of mean value of $|S_k(\alpha)|^4$:

\begin{lem}[\cite{GLZ}, Lemma 3]\label{rs}
Let $\varepsilon>0$ fixed, $k>1$, $\gamma>0$ and let $\mathcal{A}(X^{1/k};k;\gamma)$ denote the number of solutions of the inequality
\begin{align*}
|n_1^k+n_2^k-n_3^k-n_4^k|<\gamma,\qquad X^{1/k}\le n_1,n_2,n_3,n_4\le 2X^{1/k}.
\end{align*}

Then
\begin{align*}
\mathcal{A}(X^{1/k};k;\gamma)\ll\left(\gamma X^{4/k-1}+X^{2/k}\right)X^{\varepsilon}
\end{align*}
\end{lem}
\proof
See \cite{robert-sargos}, Theorem 2 with $M=X^{\frac1k}$, $\alpha=k$ and $\gamma=\delta M^k$.
\endproof

\begin{lem}[\cite{GLZ}, Lemma 4]\label{Sk_quarta_tau}
Let $k>1$, $\tau>0$. We have
\begin{align*}
\int_{-\tau}^{\tau}|S_k(\alpha)|^4\dd\alpha\ll\left(\tau X^{2/k}+X^{4/k-1}\right)X^{\varepsilon}\ll\max(\tau X^{2/k+\varepsilon},X^{4/k-1+\varepsilon}).
\end{align*}
\end{lem}

\begin{lem}[\cite{GLZ}, Lemma 10]\label{Sk_quarta}
\begin{align*}
\int_{m}|S_k(\lambda\alpha)|^4K_{\eta}(\alpha)\dd\alpha\ll\eta X^{\varepsilon}\cdot\max( X^{2/k},X^{4/k-1}).
\end{align*}
\end{lem}
%

Finally, we will use the following Lemma.
\begin{lem}\label{pnt}
\begin{align*}
&\int_m |S_1(\alpha)|^2 K_{\eta}(\alpha)\dd\alpha\ll \eta X\log X \\
&\int_m |S_2(\alpha)|^4 K_{\eta}(\alpha)\dd\alpha\ll \eta X\log^2 X.
\end{align*}
\end{lem}

\proof
The first statement comes directly from Prime Number Theorem, while the second estimation is based on Satz 3 of \cite{rieger}, p. 94.
\endproof

\section{The major arc}

Let us start from the major arc and the computation of the main term. We replace all $S_k$ defined in \eqref{S_k} with the corresponding $T_k$ defined in \eqref{T_k}. This replacement brings up some errors that we must estimate by means of Lemma \ref{s-u}, the Cauchy-Schwarz and the H\"older inequalities. We write

\begin{align*}
\mathcal{I}(\eta,\omega,\mathcal{M})=&\int_{\mathcal{M}}S_1(\lambda_1\alpha)S_2(\lambda_2\alpha)S_2(\lambda_3\alpha)S_k(\lambda_4\alpha)K_{\eta}(\alpha)e(-\omega\alpha)\dd\alpha\nonumber\\
=&\int_{\mathcal{M}}T_1(\lambda_1\alpha)T_2(\lambda_2\alpha)T_2(\lambda_3\alpha)T_k(\lambda_4\alpha)K_{\eta}(\alpha)e(-\omega\alpha)\dd\alpha\nonumber\\
&+\int_{\mathcal{M}}(S_1(\lambda_1\alpha)-T_1(\lambda_1\alpha))T_2(\lambda_2\alpha)T_2(\lambda_3\alpha)T_k(\lambda_4\alpha)K_{\eta}(\alpha)e(-\omega\alpha)\dd\alpha\nonumber\\
&+\int_{\mathcal{M}}S_1(\lambda_1\alpha)(S_2(\lambda_2\alpha)-T_2(\lambda_2\alpha))T_2(\lambda_3\alpha)T_k(\lambda_4\alpha)K_{\eta}(\alpha)e(-\omega\alpha)\dd\alpha\nonumber\\
&+\int_{\mathcal{M}}S_1(\lambda_1\alpha)S_2(\lambda_2\alpha)(S_2(\lambda_3\alpha)-T_2(\lambda_3\alpha))T_k(\lambda_4\alpha)K_{\eta}(\alpha)e(-\omega\alpha)\dd\alpha\nonumber\\
&+\int_{\mathcal{M}}S_1(\lambda_1\alpha)S_2(\lambda_2\alpha)S_2(\lambda_3\alpha)(S_k(\lambda_4\alpha)-T_k(\lambda_4\alpha))K_{\eta}(\alpha)e(-\omega\alpha)\dd\alpha\nonumber\\
=&J_1+J_2+J_3+J_4+J_5,
\end{align*}
say. For brevity, since the computations for $J_2$ and $J_3$ are similar to, but simpler than, the corresponding ones for $J_4$ and $J_5$, we will leave them to the reader.

\subsection{Main Term: lower bound for $J_1$}

As the reader might expect the main term is given by the summand $J_1$.

Let $H(\alpha)=T_1(\lambda_1\alpha)T_2(\lambda_2\alpha)T_2(\lambda_3\alpha)T_k(\lambda_4\alpha)K_{\eta}(\alpha)e(-\omega\alpha)$ so that

\begin{align*}
J_1=\int_{\z{R}}H(\alpha)\dd\alpha+\mathcal{O}\left(\int_{P/X}^{+\infty}|H(\alpha)|\dd\alpha\right).
\end{align*}

Using inequalities \eqref{k_eta} and \eqref{stima_tk},

\begin{align*}
\int_{P/X}^{+\infty}|H(\alpha)|\dd\alpha\ll& X^{-1} X^{\frac1k-1}\eta^2\int_{P/X}^{+\infty}\frac{\dd\alpha}{\alpha^4}\nonumber\ll X^{\frac1k+1}\eta^2P^{-3}=o\left(X^{\frac1k+1}\eta^2\right)
\end{align*}
provided that $P\rightarrow+\infty$. Let $D=[\delta X,X]\times[(\delta X)^{\frac12},X^{\frac12}]^2\times[(\delta X)^{\frac1k},X^{\frac1k}]$; we have
\begin{align*}
\int_{\z{R}}&H(\alpha)\dd\alpha\\&=\idotsint_D\int_{\z{R}}e((\lambda_1t_1+\lambda_2t_2^2+\lambda_3t_3^2+\lambda_4t_4^k-\omega)\alpha)K_{\eta}(\alpha)\dd\alpha\,\dd t_1\dd t_2\dd t_3\dd t_4\\&=\idotsint_D\max(0,\eta-|\lambda_1t_1+\lambda_2t_2^2+\lambda_3t_3^2+\lambda_4t_4^k-\omega)|)\dd t_1\dd t_2\dd t_3\dd t_4.
\end{align*}

Apart from trivial changes of sign, there are essentially three cases as in \cite{languasco-zaccagnini-quaternary}:
\begin{enumerate}
\item $\lambda_1>0$, $\lambda_2>0$, $\lambda_3>0$, $\lambda_4<0$ 
\item $\lambda_1>0$, $\lambda_2>0$, $\lambda_3<0$, $\lambda_4<0$ 
\item $\lambda_1>0$, $\lambda_2<0$, $\lambda_3<0$, $\lambda_4<0.$ 
\end{enumerate}

We deal with the second case, the other ones being similar: let us perform the following change of variables: $u_1=t_1-\frac{\omega}{\lambda_1}$, $u_2=t_2^2$, $u_3=t_3^2$, $u_4=t_4^k$, so that the set $D$ becomes essentially $[\delta X,X]^4$. Let us define $D'=[\delta X,(1-\delta)X]^4$ for large $X$, as a subset of $D$. The Jacobian determinant of the change of variables above is $\displaystyle\frac{1}{4k}u_2^{-\frac12}u_3^{-\frac12}u_4^{\frac1k-1}$. Then
\begin{align*}
J_1\gg&\int_{\z{R}}H(\alpha)\dd\alpha\nonumber\\
=&\idotsint_{D'}\max(0,\eta-|\lambda_1u_1+\lambda_2u_2+\lambda_3u_3+\lambda_4u_4)|)\frac{\dd u_1\dd u_2\dd u_3\dd u_4}{u_2^{\frac12}u_3^{\frac12}u_4^{1-\frac1k}}\\
\gg&X^{\frac1k-2}\hspace{-.2cm}\idotsint_{D'}\hspace{-.2cm}\max(0,\eta-|\lambda_1u_1+\lambda_2u_2+\lambda_3u_3+\lambda_4u_4)|)\dd u_1\dd u_2\dd u_3\dd u_4.
\end{align*}

Now, for $j=1,2,3$ let $a_j=\dfrac{4|\lambda_4|}{|\lambda_j|}$, $b_j=\dfrac32 a_j$ and $\mathcal{D}_j=[a_jX,b_jX]$; if $u_j\in\mathcal{D}_j$ for $j=1,2,3$ then
\begin{align*}
\lambda_1u_1+\lambda_2u_2+\lambda_3u_3\in\left[2|\lambda_4|\delta X,8|\lambda_4|\delta X\right]
\end{align*}
so that, for every choice of $(u_1,u_2,u_3)$ the interval 
\begin{align*}
[a,b]\hspace{-.1cm}=\hspace{-.1cm}\left[\frac{1}{|\lambda_4|}\left(-\eta+(\lambda_1u_1+\lambda_2u_2+\lambda_3u_3)\right),\frac{1}{|\lambda_4|}\left(\eta+(\lambda_1u_1+\lambda_2u_2+\lambda_3u_3)\right)\right]
\end{align*}
is contained in $[\delta X,(1-\delta)X]$. In other words, for $u_4\in[a,b]$ the values of $\lambda_1u_1+\lambda_2u_2+\lambda_3u_3+\lambda_4u_4$ cover the whole interval $[-\eta,\eta]$. Hence for any $(u_1,u_2,u_3)\in\mathcal{D}_1\times\mathcal{D}_2\times\mathcal{D}_3$ we have
\begin{align*}
&\int_{\delta X}^{(1-\delta)X}\max(0,\eta-|\lambda_1u_1+\lambda_2u_2+\lambda_3u_3+\lambda_4u_4|)\dd u_4\\&=|\lambda_4|^{-1}\int_{-\eta}^{\eta}\max(0,\eta-|u|)\dd u\gg\eta^2.
\end{align*}

Finally,
\begin{align*}
J_1\gg\eta^2X^{\frac1k-2}\iiint_{\mathcal{D}_1\times\mathcal{D}_2\times\mathcal{D}_3}\dd u_1\dd u_2\dd u_3\gg\eta^2X^{\frac1k-2}X^3=\eta^2X^{\frac1k+1},
\end{align*}
which is the expected lower bound.

\subsection{Bound for $J_4$}

The computations on $J_2$ and $J_3$ are similar to and simpler than the corresponding one on $J_4$, so we will skip it. Using the triangle inequality,
\begin{align*}
J_4=&\int_{\mathcal{M}}S_1(\lambda_1\alpha)S_2(\lambda_2\alpha)(S_2(\lambda_3\alpha)-T_2(\lambda_3\alpha))T_k(\lambda_4\alpha)K_{\eta}(\alpha)e(-\omega\alpha)\dd\alpha\nonumber\\
\ll&\eta^2\int_{\mathcal{M}}|S_1(\lambda_1\alpha)||S_2(\lambda_2\alpha)||S_2(\lambda_3\alpha)-T_2(\lambda_3\alpha)||T_k(\lambda_4\alpha)|\dd\alpha\nonumber\\ 
\le&\eta^2\int_{\mathcal{M}}|S_1(\lambda_1\alpha)||S_2(\lambda_2\alpha)||S_2(\lambda_3\alpha)-U_2(\lambda_3\alpha)||T_k(\lambda_4\alpha)|\dd\alpha\nonumber\\
&+\eta^2\int_{\mathcal{M}}|S_1(\lambda_1\alpha)||S_2(\lambda_2\alpha)||U_2(\lambda_3\alpha)-T_2(\lambda_3\alpha)||T_k(\lambda_4\alpha)|\dd\alpha\nonumber\\
=&\eta^2(A_4+B_4),
\end{align*}
say, where $U_2(\lambda_3\alpha)$ is given by \eqref{U_k}.

Using the trivial inequalities $|S_2(\alpha)|\ll X^{\frac12}$, \eqref{stima_tk} and then the Cauchy-Schwarz inequality,
\begin{align*}
A_4\ll &X^{\frac12}X^{\frac1k}\int_{\mathcal{M}}|S_1(\lambda_1\alpha)||S_2(\lambda_3\alpha)-U_2(\lambda_3\alpha)|\dd\alpha\nonumber\\\ll & X^{\frac12}X^{\frac1k}\left(\int_{\mathcal{M}}|S_1(\lambda_1\alpha)|^2\dd\alpha\right)^{\frac12}\left(\int_{\mathcal{M}}|S_2(\lambda_3\alpha)-U_2(\lambda_3\alpha)|^2\dd\alpha\right)^{\frac12}.
\end{align*}

Using $\int_{\mathcal{M}}|S_1(\alpha)|^2\dd\alpha\ll X\log X$ and $\int_{\mathcal{M}}|S_2(\alpha)-U_2(\alpha)|^2\dd\alpha\ll (\log X)^{-A}$ for any fixed $A$ (Lemmas \ref{s-u} and \ref{jk_stima}), we have
\begin{align*}
A_4\ll X^{\frac12+\frac1k}(X\log X)^{\frac12}(\log X)^{-\frac{A}{2}}=X^{1+\frac1k}(\log X)^{\frac12-\frac{A}{2}}=o\left(X^{\frac1k+1}\right)
\end{align*}
as long as $A>1$. Again using \eqref{stima_tk} and \eqref{t-u},
\begin{align*}
B_4\ll& X^{\frac1k}\int_{\mathcal{M}}|S_1(\lambda_1\alpha)||S_2(\lambda_2\alpha)||U_2(\lambda_3\alpha)-T_2(\lambda_3\alpha)|\dd\alpha\nonumber\\
\ll &X^{\frac1k}\int_0^{\frac1X}|S_1(\lambda_1\alpha)||S_2(\lambda_2\alpha)|\dd\alpha+X^{\frac1k+1}\int_{\frac1X}^{P/X}\alpha|S_1(\lambda_1\alpha)||S_2(\lambda_2\alpha)|\dd\alpha.
\end{align*}

Remembering that $|\alpha|\le\frac{P}{X}$ on $\mathcal{M}$ and using the H\"older inequality, trivial bounds and 
Lemma \ref{pnt}, we have
\begin{align*}
B_4\ll &X\,X^{\frac12}X^{\frac1k}\frac1X+X^{\frac1k+1}\left(\int_{1/X}^{P/X}|S_1(\alpha)|^2\dd\alpha\right)^{\frac12}\left(\int_{1/X}^{P/X}\alpha^4\dd\alpha\right)^{\frac14}\left(\int_{1/X}^{P/X}|S_2(\alpha)|^4\dd\alpha\right)^{\frac14}\\
\ll & X^{\frac12+\frac1k}+X^{\frac32+\frac1k}(\log X)^{\frac12}\left(\frac{P}{X}\right)^{\frac54}X^{\frac14}(\log X)^{\frac12}=X^{\frac12+\frac1k}P^{\frac54}\log X.
\end{align*}

Since we must have $P^{\frac54}=o(X^{\frac12}\log X)$, it follows that 
\begin{align}\label{cond2_major}
P\le X^{\frac25-\varepsilon}
\end{align}
is sufficient for our purpose.

\subsection{Bound for $J_5$}

In order to provide an estimation for $J_5$, we use \eqref{k_eta},
\begin{align*}
J_5\ll&\eta^2\int_{\mathcal{M}}|S_1(\lambda_1\alpha)||S_2(\lambda_2\alpha)||S_2(\lambda_3\alpha)||S_k(\lambda_4\alpha)-T_k(\lambda_4\alpha)|\dd\alpha
\end{align*}
and then the arithmetic-geometric inequality ($ab\le a^2+b^2$):
\begin{align*}
J_5\ll&\eta^2\sum_{j=2}^3\left(\int_{\mathcal{M}}|S_1(\lambda_1\alpha)||S_2(\lambda_j\alpha)|^2|S_k(\lambda_4\alpha)-T_k(\lambda_4\alpha)|\dd\alpha\right)
\end{align*}

The two terms are equivalent; then we consider only one of them
\begin{align*}
J_5\ll&\eta^2\int_{\mathcal{M}}|S_1(\lambda_1\alpha)||S_2(\lambda_2\alpha)|^2|S_k(\lambda_4\alpha)-T_k(\lambda_4\alpha)|\dd\alpha\nonumber\\
\ll &\eta^2\int_{\mathcal{M}}|S_1(\lambda_1\alpha)||S_2(\lambda_2\alpha)|^2|S_k(\lambda_4\alpha)-U_k(\lambda_4\alpha)|\dd\alpha\nonumber\\ 
&+\eta^2\int_{\mathcal{M}}|S_1(\lambda_1\alpha)||S_2(\lambda_2\alpha)|^2|U_k(\lambda_4\alpha)-T_k(\lambda_4\alpha)|\dd\alpha=\eta^2(A_5+B_5),
\end{align*}
say. Using trivial estimates,
\begin{align*}
A_5\ll X\int_{\mathcal{M}}|S_2(\lambda_2\alpha)|^2|S_k(\lambda_4\alpha)-U_k(\lambda_4\alpha)|\dd\alpha
\end{align*}
then using the H\"older inequality, for any fixed $A>2$ by Lemmas \ref{s-u} and \ref{jk_stima} we have
\begin{align*}
A_5\ll &X\left(\int_{\mathcal{M}}|S_2(\lambda_2\alpha)|^4\dd\alpha\right)^{\frac12}\left(\int_{\mathcal{M}}|S_k(\lambda_4\alpha)-U_k(\lambda_4\alpha)|^2\dd\alpha\right)^{\frac12}\nonumber\\\ll& X\, X^{\frac12}(\log X) \frac{P}{X}\mathcal{J}_k\left(X,\frac{X}{P}\right)^{\frac12}\ll_A X^{1+\frac1k}(\log X)^{1-\frac{A}{2}}=o\left(X^{\frac1k+1}\right),
\end{align*}
provided that $\frac{X}{P}\ge X^{1-\frac{5}{6k}+\varepsilon}$ (condition of Lemma \ref{jk_stima}), that is, 
\begin{align}\label{cond3_major}
P\le X^{\frac{5}{6k}-\varepsilon}.
\end{align}

Now we turn to $B_5$: by \eqref{t-u} we have
\begin{align*}
B_5\ll \int_0^{1/X}|S_1(\lambda_1\alpha)||S_2(\lambda_2\alpha)|^2\dd\alpha+X\int_{1/X}^{P/X}\alpha|S_1(\lambda_1\alpha)||S_2(\lambda_2\alpha)|^2\dd\alpha.
\end{align*}

Using trivial estimates and Lemma \ref{pnt}
\begin{align*}
B_5\ll&\left(\int_0^{1/X}|S_1(\lambda_1\alpha)|^2\dd\alpha\right)^{\frac12}\left(\int_0^{1/X}|S_2(\lambda_2\alpha)|^4\dd\alpha\right)^{\frac12}\nonumber\\ 
&+X\frac{P}{X}\left(\int_{1/X}^{P/X}|S_1(\lambda_1\alpha)|^2\dd\alpha\cdot\int_{1/X}^{P/X}|S_2(\lambda_j\alpha)|^4\dd\alpha\right)^{\frac12}\nonumber\\
\ll& (X\log X\cdot X\log^2 X)^{\frac12}+P(X\log X\cdot X\log^2 X)^{\frac12}\\=&X(\log X)^{\frac32}+PX(\log X)^{\frac32}.
\end{align*}

Then we need
\begin{align}\label{cond4_major}
P=o\left(X^{\frac1k-\varepsilon}\right).
\end{align}

Collecting all the bounds for $P$, that is, \eqref{cond2_major}, \eqref{cond3_major}, \eqref{cond4_major} we can take
\begin{align}\label{cond_major}
P\le \min\left(X^{\frac25-\varepsilon}, X^{\frac{5}{6k}-\varepsilon}\right).
\end{align}

In fact, if we consider \eqref{cond2_major} and \eqref{cond3_major}, we should choose the most restrictive condition between the two: if $k\le \frac{25}{12}$, $P=X^{\frac25-\varepsilon}$, otherwise, if $\frac{25}{12}<k<\frac{14}5$, $P=X^{\frac{5}{6k}-\varepsilon}$.

\section{The trivial arc}

By the arithmetic-geometric mean inequality and the trivial bound for $S_k(\lambda_4\alpha)$, we see that
\begin{align*}
|\mathcal{I}(\eta,\omega,&t)|\ll\int_{R}^{+\infty}|S_1(\lambda_1\alpha)S_2(\lambda_2\alpha)S_2(\lambda_3\alpha)S_k(\lambda_4\alpha)K_{\eta}(\alpha)|\dd\alpha\nonumber\\\ll&X^{\frac1k}\sum_{j=2}^3\int_{R}^{+\infty}|S_1(\lambda_1\alpha)|\,|S_2(\lambda_j\alpha)|^2K_{\eta}(\alpha)\dd\alpha\nonumber\\\ll& X^{\frac1k}\left(\int_{R}^{+\infty}|S_1(\lambda_1\alpha)|^2K_{\eta}(\alpha)\dd\alpha\right)^{\frac12}\left(\int_{R}^{+\infty}|S_2(\lambda_2\alpha)|^4K_{\eta}(\alpha)\dd\alpha\right)^{\frac12}\nonumber\\\ll& X^{\frac1k}\left(\int_{R}^{+\infty}\frac{|S_1(\lambda_1\alpha)|^2}{\alpha^2}\dd\alpha\right)^{\frac12}\left(\int_{R}^{+\infty}\frac{|S_2(\lambda_2\alpha)|^4}{\alpha^2}\dd\alpha\right)^{\frac12}=X^{\frac1k}C_1^{\frac12}C_2^{\frac12},
\end{align*}
say. Using the PNT and the periodicity of $S_1(\alpha)$, we have
\begin{align}\label{c1}
C_1=&\int_{R}^{+\infty}\frac{|S_1(\lambda_1\alpha)|^2}{\alpha^2}\dd\alpha\ll\int_{\lambda_1R}^{+\infty}\frac{|S_1(\alpha)|^2}{\alpha^2}\dd\alpha\nonumber\\\ll&\sum_{n\ge|\lambda_1| R}\frac{1}{(n-1)^2}\int_{n-1}^n |S_1(\alpha)|^2\dd\alpha\ll\frac{X\log X}{|\lambda_1| R}.
\end{align}

Now using Lemma \ref{pnt},
\begin{align}\label{c2}
C_2=&\int_{R}^{+\infty}\frac{|S_2(\lambda_2\alpha)|^4}{\alpha^2}\dd\alpha\ll\int_{\lambda_1R}^{+\infty}\frac{|S_2(\alpha)|^4}{\alpha^2}\dd\alpha\nonumber\\\ll&\sum_{n\ge|\lambda_1| R}\frac{1}{(n-1)^2}\int_{n-1}^n |S_2(\alpha)|^4\dd\alpha\ll\frac{X\log^2 X}{|\lambda_1| R}.
\end{align}

Collecting \eqref{c1} and \eqref{c2},
\begin{align*}
|\mathcal{I}(\eta,\omega,t)|\ll X^{\frac1k}\left(\frac{X\log X}{R}\right)^{\frac12}\left(\frac{X\log^2 X}{R}\right)^{\frac12}\ll\frac{X^{1+\frac1k}(\log X)^{\frac32}}{R}.
\end{align*}

Hence, remembering that $|\mathcal{I}(\eta,\omega,t)|$ must be $o\left(\eta^2X^{\frac1k+1}\right)$, i.e. of the main term, the choice
\begin{align}\label{cond_trivial}
R=\frac{\log^2 X}{\eta^2}
\end{align}
is admissible.

\section{The minor arc}

In \cite{languasco-zaccagnini-ternary} Lemma 3 it is proven that the measure of the set where $|S_1(\lambda_1\alpha)|^{\frac12}$ and $|S_2(\lambda_2\alpha)|$ are both large for $\alpha\in m$ is small, exploiting the fact that the ratio $\lambda_1/\lambda_2$ is irrational.

\begin{lem}[Vaughan \cite{vaughan1997hardy}, Theorem 3.1]\label{vaughan}
Let $\alpha$ be a real number and $a,q$ be positive integers satisfying $(a,q)=1$ and $\left|\alpha-\frac{a}{q}\right|<\frac{1}{q^2}$. Then
\begin{align*}
S_1(\alpha)\ll\left(\frac{X}{\sqrt{q}}+\sqrt{X q}+X^{\frac45}\right)\log^4 X.
\end{align*}
\end{lem}

We now state some considerations about Lemmas \ref{vaughan}:

\begin{cor}[Liu-Sun \cite{liu2013diophantine}, Corollary 2.7]\label{harman_s1}
Suppose that $X\ge Z\ge X^{1-\frac15+\varepsilon}$ and $|S_1(\lambda_1\alpha)|>Z$. Then there are coprime integers $(a,q)=1$ satisfying
\begin{align*}
1\le q\le \left(\frac{X^{1+\varepsilon}}{Z}\right)^2,\qquad |q\lambda_1\alpha-a|\ll\left(\frac{X^{\frac12+\varepsilon}}{Z}\right)^2.
\end{align*}
\end{cor}

\begin{lem}[Wang-Yao \cite{wang-yao}, Lemma 1]\label{harman_s2}
Suppose that $X^{\frac12}\ge Z\ge X^{\frac12-\frac{1}{14}+\varepsilon}$ and $|S_2(\lambda_2\alpha)|>Z$. Then there are coprime integers $(a,q)=1$ satisfying
\begin{align*}
1\le q\le \left(\frac{X^{\frac12+\varepsilon}}{Z}\right)^4,\qquad |q\lambda_2\alpha-a|\ll X^{-1}\left(\frac{X^{\frac12+\varepsilon}}{Z}\right)^4
\end{align*}
\end{lem}

Let us now split $m$ into two subsets $\tilde{m}$ and $m^*=m\backslash\tilde{m}$. In turn $\tilde{m}=m_1\cup m_2$, where
\begin{align*}
&m_1=\{\alpha\in m\,:\ |S_1(\lambda_1\alpha)|\le X^{1-\frac17+\varepsilon}\} \\
&m_2=\{\alpha\in m\,:\ |S_2(\lambda_2\alpha)|\le X^{\frac12-\frac1{14}+\varepsilon}\}. 
\end{align*}

Using the H\"older inequality, Lemma \ref{Sk_quarta} and the definition of $m_1$ we obtain
\begin{align}\label{cond_m1}
|\mathcal{I}(\eta,\omega,&m_1)|\ll\int_{m_1}|S_1(\lambda_1\alpha)||S_2(\lambda_2\alpha)||S_2(\lambda_3\alpha)||S_k(\lambda_4\alpha)|K_{\eta}(\alpha)\dd\alpha\nonumber\\ 
\ll&\left(\max_{\alpha\in m_1}{|S_1(\lambda_1\alpha)|}\right)^{\frac12}\left(\int_{m_1}|S_1(\lambda_1\alpha)|^2K_{\eta}(\alpha)\dd\alpha)\right)^{1/4}\nonumber\\
&\prod_{i=2}^3\left(\int_{m_1}|S_2(\lambda_i\alpha)|^4K_{\eta}(\alpha)\dd\alpha)\right)^{1/4}\left(\int_{m_1}|S_k(\lambda_4\alpha)|^4K_{\eta}(\alpha)\dd\alpha)\right)^{1/4}\nonumber\\
\ll& X^{\frac{3}{7}+\varepsilon}(\eta X\log X)^{\frac14}(\eta X\log^2 X)^{\frac12}\left(\eta X^{\varepsilon}\max(X^{\frac2k},X^{\frac4k-1})\right)^{\frac14}\nonumber\\
=&\eta X^{\frac{33}{28}+2\varepsilon}\max(X^{\frac1{2k}},X^{\frac1k-\frac14}).
\end{align}

Using the H\"older inequality, Lemma \ref{Sk_quarta} and the definition of $m_2$ we obtain
\begin{align}\label{cond_m2}
|\mathcal{I}(\eta,\omega,&m_2)|\ll\int_{m_2}|S_1(\lambda_1\alpha)||S_2(\lambda_2\alpha)||S_2(\lambda_3\alpha)||S_k(\lambda_4\alpha)|K_{\eta}(\alpha)\dd\alpha\nonumber\\ 
\ll&\max_{\alpha\in m_2}{|S_2(\lambda_2\alpha)|}\left(\int_{m_2}|S_1(\lambda_1\alpha)|^2K_{\eta}(\alpha)\dd\alpha)\right)^{1/2}\nonumber\\
&\left(\int_{m_2}|S_2(\lambda_3\alpha)|^4K_{\eta}(\alpha)\dd\alpha)\right)^{1/4}\left(\int_{m_2}|S_k(\lambda_4\alpha)|^4K_{\eta}(\alpha)\dd\alpha)\right)^{1/4}\nonumber
\\\ll& X^{\frac{3}{7}+\varepsilon}(\eta X\log X)^{\frac12}(\eta X\log^2 X)^{\frac14}\left(\eta X^{\varepsilon}\max(X^{\frac2k},X^{\frac4k-1})\right)^{\frac14}\nonumber\\
=&\eta X^{\frac{33}{28}+2\varepsilon}\max(X^{\frac1{2k}},X^{\frac1k-\frac14}).
\end{align}

Both \eqref{cond_m1} and \eqref{cond_m2} must be $o\left(\eta^2X^{1+\frac1k}\right)$, consequently it is clear that for $1<k<2$, $\eta$ is a negative power of $X$ independently from the value of $k$. The we have the following most restrictive condition for $k\ge 2$:
\begin{align*}
\eta=\infty\left(X^{{-\frac{15-5k}{28k}}+\varepsilon}\right).
\end{align*}

It remains to discuss the set $m^*$ in which the following bounds hold simultaneously
\begin{align*}
|S_1(\lambda_1\alpha)|>X^{\frac67+\varepsilon},\, |S_2(\lambda_2\alpha)|>X^{\frac37+\varepsilon},\, \frac{P}{X}=X^{-\frac35}<|\alpha|\le \frac{\log^2 X}{\eta^2}=R.
\end{align*} 

Following the dyadic dissection argument as in \cite{harman2004ternary} we divide $m^*$ into disjoint sets $E(Z_1,Z_2,y)$ in which, for $\alpha\in E(Z_1,Z_2,y)$, we have
\begin{align*}
Z_1<|S_1(\lambda_1\alpha)|\le 2Z_1,\qquad Z_2<|S_2(\lambda_2\alpha)|\le 2Z_2,\qquad y<|\alpha|\le 2y
\end{align*}
where $Z_1=2^{k_1}X^{\frac78+\varepsilon}$, $Z_2=2^{k_2}X^{\frac7{16}+\varepsilon}$ and $y=2^{k_3}X^{-\frac34-\varepsilon}$ for some non-negative integers $k_1,k_2,k_3$.

It follows that the disjoint sets are, at the most, $\ll\log^3 X$. Let us define $\mathcal{A}$ as a shorthand for the set $E(Z_1,Z_2,y)$; we have the following result about the Lebesgue measure of $\mathcal{A}$ following the same lines of Lemma 6 in \cite{mu2016diophantine}:
\begin{lem}\label{lemma_mu}
We have $\mu(\mathcal{A})\ll yX^{\frac{18}7+6\varepsilon}Z_1^{-2}Z_2^{-4}$,
where $\mu(\cdot)$ denotes the Lebesgue measure.
\end{lem}
\proof
If $\alpha\in\mathcal{A}$, by Corollaries \ref{harman_s1} and Lemma \ref{harman_s2} there are coprime integers $(a_1,q_1)$ and $(a_2,q_2)$ such that
\begin{align}\label{misura1}
&1\le q_1\ll\left(\frac{X^{1+\varepsilon}}{Z_1}\right)^2,\qquad |q_1\lambda_1\alpha-a_1|\ll\left(\frac{X^{\frac12+\varepsilon}}{Z_1}\right)^2\nonumber\\&1\le q_2\ll\left(\frac{X^{\frac12+\varepsilon}}{Z_2}\right)^4,\qquad |q_2\lambda_2\alpha-a_2|\ll X^{-1}\left(\frac{X^{\frac12+\varepsilon}}{Z_2}\right)^4.
\end{align}

We remark that $a_1a_2\neq0$ otherwise we would have $\alpha\in\mathcal{M}$. Recalling the definitions of $Z_1$ and $Z_2$:
\begin{align}\label{alpha_esplicito}
&q_1^{-1}\left(\frac{X^{\frac12+\varepsilon}}{Z_1}\right)^2\ll\alpha\ll q_1^{-1}\left(\frac{X^{\frac12+\varepsilon}}{Z_1}\right)^2\nonumber\\&q_2^{-1}X^{-1}\left(\frac{X^{\frac12+\varepsilon}}{Z_2}\right)^4\ll\alpha\ll q_2^{-1}X^{-1}\left(\frac{X^{\frac12+\varepsilon}}{Z_2}\right)^4;
\end{align}
then, if $q_1=q_2=1$ we get:
\begin{align*}
\alpha\gg\frac{X}{X^{\frac{12}7+2\varepsilon}}=X^{-\frac57+2\varepsilon}.
\end{align*}

Note that, by the conditions on $P$ (see \eqref{cond_major}), there is no gap between the major arc and the minor arc. Now, we can further split $m^*$ into sets $I(Z_1,Z_2,y,Q_1,Q_2)$ where, on each set, \mbox{$Q_j\le q_j\le 2Q_j$}. If we explicit $\alpha$ as in \eqref{alpha_esplicito}, we obtain the following inequalities:
\begin{align*}
&Q_1^{-1}\left(\frac{X^{\frac12+\varepsilon}}{Z_1}\right)^2\ll\alpha\ll Q_1^{-1}\left(\frac{X^{\frac12+\varepsilon}}{Z_1}\right)^2\nonumber\\&Q_2^{-1}X^{-1}\left(\frac{X^{\frac12+\varepsilon}}{Z_2}\right)^4\ll\alpha\ll Q_2^{-1}X^{-1}\left(\frac{X^{\frac12+\varepsilon}}{Z_2}\right)^4.
\end{align*}

As the inequalities \eqref{misura1} hold simultaneously, the measure of $I$ can be bounded with the minimum of the two:
\begin{align*}
\mu(I)\ll\min\left(Q_1^{-1}\left(\frac{X^{\frac12+\varepsilon}}{Z_1}\right)^2,Q_2^{-1}X^{-1}\left(\frac{X^{\frac12+\varepsilon}}{Z_2}\right)^4\right).
\end{align*}

Taking the geometric mean ($\min(a,b)\le\sqrt{a}\sqrt{b}$) we can write
\begin{align}\label{lower_Q}
\mu(I)&\ll Q_1^{-\frac12}Q_2^{-\frac12}X^{-\frac12}\left(\frac{X^{\frac12+\varepsilon}}{Z_1}\right)\left(\frac{X^{\frac12+\varepsilon}}{Z_2}\right)^2\ll \frac{X^{1+3\varepsilon}}{Q_1^{\frac12}Q_2^{\frac12}Z_1Z_2^2}.
\end{align}

Now we need a lower bound for $Q_1^{\frac12}Q_2^{\frac12}$: by \eqref{misura1}
\begin{align*}
\left|a_2q_1\frac{\lambda_1}{\lambda_2}-a_1q_2\right|&=\left|\frac{a_2}{\lambda_2\alpha}(q_1\lambda_1\alpha-a_1)-\frac{a_1}{\lambda_2\alpha}(q_2\lambda_2\alpha-a_2)\right|\nonumber\\&\ll q_2|q_1\lambda_1\alpha-a_1|+q_1|q_2\lambda_2\alpha-a_2|\ll\nonumber
\\&\quad Q_2\left(\frac{X^{\frac12+\varepsilon}}{Z_1}\right)^2+Q_1X^{-1}\left(\frac{X^{\frac12+\varepsilon}}{Z_2}\right)^4.
\end{align*}

Remembering that $Q_1\ll \left(\frac{X^{1+\varepsilon}}{Z_1}\right)^2$, $Q_2\ll\left(\frac{X^{\frac12+\varepsilon}}{Z_2}\right)^4$, $Z_1\gg X^{\frac67+\varepsilon}$, $Z_2\gg X^{\frac37+\varepsilon}$, 

\begin{align}\label{harman_stima}
\left|a_2q_1\frac{\lambda_1}{\lambda_2}-a_1q_2\right|&\ll\left(\frac{X^{\frac12+\varepsilon}}{X^{\frac37+\varepsilon}}\right)^4\left(\frac{X^{\frac12+\varepsilon}}{X^{\frac67+\varepsilon}}\right)^2+\left(\frac{X^{1+\varepsilon}}{X^{\frac67+\varepsilon}}\right)^2X^{-1}\left(\frac{X^{\frac12+\varepsilon}}{X^{\frac37+\varepsilon}}\right)^4\nonumber\\&\ll\frac{X^{2+4\varepsilon}X^{1+2\varepsilon}}{X^{\frac{12}7+4\varepsilon}X^{\frac{12}7+2\varepsilon}}\ll X^{-\frac37-6\varepsilon}<\frac{1}{4q}
\end{align}
since $X=q^{7/3}$. We would like that $|a_2q_1|\ge q$ so that, recalling that $q$ is the denominator of a convergent of $\lambda_1/\lambda_2$, we could apply the Legendre's law of best approximation for continued fractions: in our case it must be
\begin{align*}
X^{-\frac37-6\varepsilon}<\frac{1}{4 q}.
\end{align*}

It turns out that for any pair $\alpha$, $\alpha'$ having distinct associated products $a_2q_1$ (see Watson \cite{watson}),
\begin{align*}
|a_2(\alpha)q_1(\alpha)-a_2(\alpha')q_1(\alpha')|\ge q;
\end{align*}
thus, by the pigeon-hole principle, there is at most one value of $a_2q_1$ in the interval $[rq,(r+1)q)$ for any positive integer $r$. $a_2q_1$ determines $a_2$ and $q_1$ to within $X^{\varepsilon}$ possibilities (from the bound for the divisor function) and consequently also $a_2q_1$ determines $a_1$ and $q_2$ to within $X^{\varepsilon}$ possibilities from \eqref{harman_stima}.

Hence we got a lower bound for $q_1q_2$, remembering that in our shorthand $Q_j\le q_j\le 2Q_j$:
\begin{align*}
q_1q_2=a_2q_1\frac{q_2}{a_2}\gg\frac{rq}{|\alpha|}\gg rqy^{-1}
\end{align*}
and finally by \eqref{lower_Q}
\begin{align*}
\mu(I)\ll X^{1+3\varepsilon}{Z_1^{-1}Z_2^{-2}}r^{-\frac12}q^{-\frac12}y^{\frac12}.
\end{align*}

Inside the interval $[rq,(r+1)q)$, $rq\le |a_2q_1|$ and, in turn from \eqref{misura1}, $a_2\ll q_2|\alpha|$, then
\begin{align*}
&rq\ll q_1q_2|\alpha|\ll\left(\frac{X^{1+\varepsilon}}{Z_1}\right)^2\left(\frac{X^{\frac12+\varepsilon}}{Z_2}\right)^4y\ll yX^{4+6\varepsilon}Z_1^{-2}Z_2^{-4}\nonumber\\&\Rightarrow r\ll q^{-1}yX^{4+6\varepsilon}Z_1^{-2}Z_2^{-4}.
\end{align*}

Now, we sum on every interval to get an upper bound for the measure of $\mathcal{A}$:
\begin{align*}
\mu(\mathcal{A})\ll X^{1+3\varepsilon}{Z_1^{-1}Z_2^{-2}}q^{-\frac12}y^{\frac12}\sum_{1\le r\ll q^{-1}yX^{4+6\varepsilon}Z_1^{-2}Z_2^{-4}}r^{-\frac12}.
\end{align*}

By partial summation on the generalized harmonic series, 
\begin{align*}
\sum_{1\le r\ll q^{-1}yX^{4+6\varepsilon}Z_1^{-2}Z_2^{-4}}r^{-\frac12}\ll (q^{-1}yX^{4+6\varepsilon}Z_1^{-2}Z_2^{-4})^{\frac12}
\end{align*}
then
\begin{align*}
\mu(\mathcal{A})&\ll yX^{3+6\varepsilon}Z_1^{-2}Z_2^{-4}q^{-1}\ll yX^{3+6\varepsilon}Z_1^{-2}Z_2^{-4}X^{-\frac37}\ll yX^{\frac{18}7+6\varepsilon}Z_1^{-2}Z_2^{-4}.
\end{align*}
\endproof

Using Lemma \ref{lemma_mu} we finally are able to get a bound for $\mathcal{I}(\eta,\omega,\mathcal{A})$:
\begin{align}\label{final}
|\mathcal{I}(\eta,\omega,\mathcal{A})|\ll&\int_{\mathcal{A}}|S_1(\lambda_1\alpha)||S_2(\lambda_2\alpha)||S_2(\lambda_3\alpha)||S_k(\lambda_4\alpha)|K_{\eta}(\alpha)\dd\alpha\nonumber\\
\ll&\left(\int_{\mathcal{A}}|S_1(\lambda_1\alpha)S_2(\lambda_2\alpha)|^2K_{\eta}(\alpha)\dd\alpha\right)^{\frac12}\left(\int_{\mathcal{A}}|S_2(\lambda_3\alpha)|^4K_{\eta}(\alpha)\dd\alpha\right)^{\frac14}\nonumber\\
&\left(\int_{\mathcal{A}}|S_k(\lambda_4\alpha)|^4K_{\eta}(\alpha)\dd\alpha\right)^{\frac14}\nonumber\\
\ll&\left(\min\left(\eta^2,\frac{1}{y^2}\right)\right)^{\frac12}\left((Z_1Z_2)^2\mu(\mathcal{A})\right)^{\frac12}\left(\eta X \log^2 X\right)^{\frac14}\left(\eta X^{\varepsilon}\max(X^{\frac2k},X^{\frac4k-1})\right)^{\frac14}\nonumber\\
\ll&\eta Z_2^{-1}X^{\frac97+2\varepsilon}X^{\frac14+2\varepsilon}\max(X^{\frac1{2k}},X^{\frac1k-\frac14})\ll \eta X^{\frac67+2\varepsilon}X^{\frac14+2\varepsilon}\max(X^{\frac1{2k}},X^{\frac1k-\frac14})
\nonumber\\
\ll&\eta X^{\frac{31}{18}+4\varepsilon}\max(X^{\frac1{2k}},X^{\frac1k-\frac14})
\end{align}
so $\eta=\infty\left(X^{{-\frac{14-5k}{28k}}+\varepsilon}\right)$ is the optimal choice.

\section*{Acknowledgement}

I wish to thank my PhD supervisor Prof. Alessandro Zaccagnini for his professional guidance and valuable support and Prof. Alessandro Languasco for his fruitful suggestions which improved my article. This is a part of my PhD dissertation.

\end{document}